\input AHTOHFIE.STY
\def\F{{\cal F}}

\def\repE{\Upsilon_{\rm e}}
\def\repESym{\Upsilon_{\rm e}^{\rm sym}}

\UDC{\hskip-0.1em
519.157.1
+519.175.1
+519.176
}

\title{
Invariant systems of weighted representatives
}
\author{%
Anton A. Klyachko
\quad
Mikhail S. Terekhov
}
\address{
\myAddressWC
\qquad\qquad
kombox.ver.2.0@yandex.ru
}

\grants{\RSF 22-11-00075}

\abstract{%
It is known that, if removing some $n$ edges from a graph $\Gamma$
destroys all subgraphs
isomorphic to a given finite graph $K$,
then
all subgraphs
isomorphic to $K$ can
be destroyed by
removing at most $|E(K)|\cdot n$ edges,
which form a set invariant with respect to all automorphisms of $\Gamma$.
We construct the first examples of (connected)
graphs $K$ for which this estimate is not
sharp.
Our arguments are based on a ``weighted analogue"
of an earlier known estimate for the cost of symmetry.
}

\s 1.
Introduction

\disp{%
Suppose that we
can choose $n$ vertices of a graph in such a way that each
100-cycle contains a chosen vertex.
How many vertices should be chosen if we want
to make the choice
\emph{fair}, i.e., we
want
the set of chosen vertices be invariant with respect to
all automorphisms of the graph\/\rm?
}%
It is known that the answer is {\sl at most $100n$ suffice}
(and this estimate is sharp). Moreover,
the following general fact was obtained in [KL21].

\Theorem KL {\rm[KL21]}.
Suppose that a group $G$ acts on a set $U$,
$\F$ is a $G$-invariant family
of finite
subsets of $U$ of
uniformly bounded
cardinality,
and $X\subseteq U$ is a
finite
system of representatives for this
family \(i.e., $X\cap F\ne\emptyset$ for any~$F\in\F$\).
Then there exists a $G$-invariant system of representatives $Y$ such that
$
|Y|\le |X|\cdot
\max
\limits_{F\in\F}|F|.
$

\noindent
This fact can be called a combinatorial analogue of the
(algebraic) Khukhro--Makarenko theorem
[KhM07a] (see also [KhM07b], [KlMe09], and [KlMi15]).
A survey of algebraic Khukhro--Makarenko type theorems
can be found in [KL21].

First, we prove the following generalisation of Theorem KL (in Section 3).

\Theorem on symmetrisation of systems of multiple representatives.
Suppose that a group $G$ acts on a set~$U$, $\F$ is a $G$-invariant
family of finite
subsets of
$U$ of uniformly bounded
cardinality,
and $X\subseteq U$ is a
finite
\emph{system of $k$-multiple
representatives} for this
family \(i.e. $|X\cap F|\ge k$ for any~$F\in\F$\).
Then there exists a $G$-invariant
system of $k$-multiple representatives $Y$ such that
$
k\cdot|Y|\le |X|\cdot
\max
\limits_{F\in\F}|F|.
$
\newline
Moreover,
as $Y$, we can take the union of all $G$-orbits $G\o u$ such that
$|(G\o u)\cap X|
\cdot
\max\limits_{F\in\F}|F|
\ge
|G\o u|\cdot k
$.
In particular,
$
k\cdot|Y|\le|X\cap Y|\cdot
\max\limits_{F\in\F}|F|
$.

For example, if a graph $\Gamma$ possesses a finite set of vertices $X$ 
containing at least two vertex from each 100-cycle, then we can choose an 
$(\Aut\Gamma)$-invariant set of vertices $Y$ containing at least two 
vertex from each 100-cycle and such that $|Y|\le50|X|$.

In Section 4, we prove
an even more general version of Theorem KL, which makes it possible, e.g.,
to solve the following ``applied" problem:
\disp{\it\narrower\narrower\hfuzz2.5pt
to form a student council, we have, from every ten students,
the first of which is respected by the nine others, to choose 
\-
either this respected student (the first one),
\-
or at least three respecting student (from the nine remaining ones).
\enditem
We want to minimise the cardinality
of the student council.
If we wish also to choose fairly, then what is the cost of fairness?
}%
There arises a directed graph with weights. In this case, this is the
ten-vertex star; the
central vertex has weight one, and the nine remaining vertices
are of weight $1\over3$.
There is also a (large) directed graph $\Gamma$ (with no
weights) describing who respects whom. We have to choose
some vertices of $\Gamma$ in such a way that each ten-vertex star is 
presented with weight at least one. If we do not care about fairness, then 
we can
choose a student council consisting of, say, forty persons.
If we also want to be fair
(i.e., we want the student council to be invariant with
respect to $\Aut\Gamma$), how many students must be chosen in the worst 
case?  The answer is $\le40\cdot\(1+9\cdot{1\over3}\)=160$, because of the 
following natural generalisation of Theorem~KL
(which is proven in Section 4).

\Theorem on symmetrisation of
systems of weighted representatives.
Suppose that a group $G$ acts on a set~$U$,
$W\subset\Q$ is a finite set \(of weights\),
and
$\F$ is a $G$-invariant family of
nonnegative \(\emph{weight}\) functions~$U\to W$
with finite supports
\(where $G$-invariance means
that, for each $F\in\F$ and $g\in G$,
the function $u\mapsto F(g\o u)$ also lies in~$\F$\).
Let $X\subseteq U$ be a
finite
\emph{system of weighted representatives} for this
family,
i.e. $\sum\limits_{x\in X}F(x)\ge1$ for any~$F\in\F$.
Then there exists a $G$-invariant system of weighted representatives
$Y\subseteq U$ such that
$
|Y|\le|X|\cdot
\max\limits_{F\in\F}\sum\limits_{u\in U}F(u).
$
\newline
Moreover,
as $Y$, we can take the union of all $G$-orbits $G\o u$ such, that
$|(G\o u)\cap X|
\cdot
\max\limits_{F\in\F}\sum\limits_{u\in U}F(u)
\ge
|G\o u|
$.
In particular,
$
|Y|\le|X\cap Y|\cdot
\max\limits_{F\in\F}\sum\limits_{u\in U}F(u)
$.

We apply this result to solve an open problem on graphs
([KL21], Question 2).
Theorem KL implies immediately the following fact.

\Corollary {\rm(of Theorem KL) [KL21]}.
Let $\Gamma$ be a graph, and let $K$ be a finite graph. Then
\item{\rm (v)}
if $\Gamma$ contains a finite set of vertices $X$
such that each subgraph of $\Gamma$
isomorphic to $K$ has at least one
vertex from~$X$, then
$\Gamma$ contains an $(\Aut\Gamma)$-invariant set of vertices $Y$
such that again
each subgraph isomorphic to $K$ have at least one
vertex from $Y$ and
$
|Y|\le |X|\cdot(\hbox{\rm the number of vertices of $K$});
$
\item{\rm (e)}
if $\Gamma$ contains a
finite set of
edges $X$ such that each subgraph isomorphic to $K$ have at least one
edge from $X$, then
$\Gamma$ contains an $(\Aut\Gamma)$-invariant
set of edges~$Y$ such that again
each subgraph isomorphic to $K$ has at least one
edge from $Y$ and
$
|Y|\le |X|\cdot(\hbox{\rm the number of edges of $K$}).
$

\disp{\sl
Are these estimates sharp for each given connected
graph $K$\?
}%
\item{(v)}
For estimate (v), the
answer is positive in the general case [KL21],
but negative if we require additionally the connectedness of
$\Gamma$ [T22] (and the minimum counterexample is a five-vertex tree
$K=D_5$).
\item{(e)}
We show that, for estimate (e), the answer is also
negative (even if the connectedness of $\Gamma$ is not required),
answering thereby an open question from [KL21].
We apply the theorem on symmetrisation of systems of weighted 
representatives (and the Dulmage--Mendelsohn decomposition [DM58]) to
construct a large class of examples of such \emph{non-edge-costly} graphs
$K$ (Fig. 1), see the following section (which contains also a formal
definition of ``sharpness").

\goodbreak
\bigskip
\centerline{\input 1.PIC}
\nobreak%
\centerline{Fig. 1}%
\goodbreak
\bigskip

The word \emph{graph} in this paper means an undirected
(possibly infinite)
graph without loops
and multiple edges.
The symbols $V(\Gamma)$ and $E(\Gamma)$ denote the sets of
vertices and edges of a graph $\Gamma$.

The authors thank
Alexander Skutin for valuable remarks;
we thank also
the Theoretical Physics and Mathematics Advancement Foundation ``BASIS".

\s 2.
Edge non-costly ``tadpoles-graphs"

The \emph{edge representativeness} $\repE(K,\Gamma)$
of a graph $K$ in a graph $\Gamma$
[KL21]
is the
minimal integer $n$ such that
$\Gamma$ contains a set of edges $X$ of cardinality~$n$
such that
$$
\hbox{\sl each subgraph of $\Gamma$ isomorphic to $K$
contains an edge from $X$.}
\eqno{(1)}
$$
The \emph{symmetric edge representativeness}
$\repESym(K,\Gamma)$ of a graph $K$ in a graph~$\Gamma$
[KL21]
is the
minimal
integer $n$ such that $\Gamma$ an $(\Aut\Gamma)$-invariant
set of edges $X$ of cardinality~$n$ satisfying $(1)$.

Clearly, $\repE(K,\Gamma)\le\repESym(K,\Gamma)$.
The corollary of Theorem KL
says that
$\repESym(K,\Gamma)\le\repE(K,\Gamma)\cdot|E(K)|$.
A graph~$K$ is called [KL21]
\emph{costly in the sense symmetric edge representativeness}
or simply
\emph{edge-costly} if
$$
\hbox{
$\forall m\in\Z$ there exists a graph
$\Gamma_m$ such that
}
\repESym(K,\Gamma_m)=
\repE(K,\Gamma_m)\cdot(\hbox{the number of edges of $K$})\ge m.
\eqno{(2)}
$$
Thus, a graph $K$ is edge-costly if the estimate from
the corollary (about edges) of
Theorem KL
cannot be improved for~$K$.

\Th on the edge representativeness of ``tadpoles-graphs".
If a finite connected graph $K$ contains precisely one vertex
of degree one,
and the graph obtained from $K$ by removing this vertex and
the incident edge is vertex-transitive\/ {\rm(Fig. 1)},
then
$
\repESym(K,\Gamma)
\le
\bigl(|E(K)|-1\bigr)
\cdot\repE(K,\Gamma)
$
for any graph $\Gamma$\;
in particular, $K$ is not edge-costly.

\Proof
Let $K_0$ be the vertex-transitive graph obtained from $K$
by removing the ``tail".
The first observation is that
\disp{\sl\hfuzz10cm
a set of representatives $X$ for
the family $\{E(\~K)\;|\;\Gamma\supseteq\~K\iso K\}$
must be simultaneously a system of weighted representatives
for the
family $\F'$ of weight functions of the following form
\newline\phantom{.}
\newline
$
\phantom{.}\qquad
F_{K'K''}\:E(\Gamma)\to\Q,
\quad
F_{K'K''}(e)=
\cases{
1,&if $e\in K'\cap K''$
\cr
1/2,&if $e\in (K'\cup K'')\setminus(K'\cap K'')$
\cr
0,&if $e\notin (K'\cup K'')$
\cr
},
$
\newline\phantom{.}
\newline
where $K'$ and $K''$ are subgraphs of $\Gamma$ isomorphic to $K_0$
such that $K'\cup K''$ is a connected graph and
$V(K')\ne V(K'\cup K'')\ne V(K'')$
\rm\(see Fig. 2, where $K$ is the tailed triangle\).
}%

\goodbreak
\bigskip
\centerline{\input 2.PIC}
\nobreak%
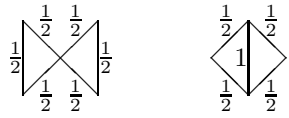
\centerline{Fig. 2}%
\goodbreak
\bigskip

Indeed, to destroy all subgraphs of $K'\cup K''$ isomorphic
to $K$, we have either to remove an edge from $K'\cap K''$
or to remove
at least two edges, because removing one edge $e$, say,
from~$K'\setminus(K'\cap K'')$ does not destroy $K''$, and a tail for $K''$
also remains intact, since
\-
the graph $K'\setminus\{e\}$ is connected
(as the edge connectedness of a 
finite
vertex-transitive graph equals
the degree of a vertex
[Ma71],
i.e., greater than one, by the condition of the theorem);
\-
therefore, the graph $K''\cup K'\setminus\{e\}$ is connected;
\-
there exists a vertex $v$ lying in $K''$, but not in $K'$.

\medskip\noindent
For each function $F\in\F'$, we have
$
\sum\limits_{u}F(u)=|E(K_0)|=|E(K)|-1,
$
hence,
the theorem on symmetrisation of systems of weighted representatives
implies that
some
symmetric system of weighted representatives $Y'$ for
$\F'$ contains at most $(|E(K)|-1)\cdot|X\cap Y'|$ edges.

Now, consider the graph $\Gamma'=\Gamma\setminus Y'$,
on which the group $G=\Aut\Gamma$ acts
(because the set $Y'$ is invariant).
This graph has the following properties:
\item{1)}
the set $X'=X\setminus(Y'\cap X)$
is a system of representatives
for the edges of
subgraphs isomorphic to $K$;
\item{2)}
the vertex-sets of any two subgraphs (of $\Gamma'$) isomorphic
to $K_0$
either coincide or does not intersect.

\enditem
By virtue of 1),
\dispno{\sl\hfuzz1cm
it remains to find a symmetric~~system of representatives $Y''$
for edges
\newline
of subgraphs isomorphic to
$K$ in $\Gamma'$
such
that $|Y''|\le(|E(K)|-1)\cdot|X'|$\;
}(*)%
because then $Y=Y'\cup Y''$ would be a required symmetric system
representatives in $\Gamma$ (as $\Gamma\setminus Y$ does not
contain subgraphs isomorphic to $K$) and
$$
|Y|=|Y'|+|Y''|\le(|E(K)|-1)\cdot|X\cap Y'|+(|E(K)|-1)\cdot|X'|=
(|E(K)|-1)\cdot|X|.
$$
\goodbreak
\noindent
For a set $S\subseteq V(\Gamma')$ of cardinality $|V(K_0)|$,
put $m_S=\Upsilon_e(K_0,\gp S)$, where $\gp S$ is the induced subgraph 
of $\Gamma'$ with the vertex set $S$, and let 
$A_S$ be an (abstract) set of cardinality $m_S$.

\noindent
Consider the following auxiliary biparite graph $\Delta$:
\-
the fraction $A$ is $\bigsqcup\limits_S A_S$
(henceforth $S$ ranges over all sets $S\subseteq V(\Gamma')$ 
of cardinality $|V(K_0)|$).
\-
the vertices of the fraction $B$ are
the edges of $\Gamma'$
lying in no subgraphs isomorphic to~$K_0$;
\-
$a\in A_S$ and $b\in B$ are joined by an edge
if, in $\Gamma'$,
precisely one endpoint of the edge $b$ lies in $S$.
\enditem
The set 
$
Q=\(\bigsqcup\{A_S\mid |X'\cap S|\ge m_S\}\)\sqcup(X'\cap B)
$
is a vertex cover of $\Delta$
(i.e., it represents all edges). 
Indeed,
if an edge with endpoints $a\in A_S$ and $b\in B$ is not represented,
then the graph $\gp S\setminus X'$ has a subgraph $G$ isomorphic to~$K_0$, 
and $G\cup b$ is a subgraph of $\Gamma'\setminus X'$ isomorphic to $K$, 
which contradicts the definition of $X'$.

Now, $|Q|=\sum|A_S|+|X'\cap B|$, 
where the sum ranges over all $S$ such that $|X'\cap S|\ge m_S=|A_S|$;
so, each term of this sum is at most $|X'\cap S|$
and $|Q|\le|X'|$. 

The
Dulmage--Mendelsohn decomposition [DM58]
(see also
[LP86],
Theorem 3.2.4)
implies that
\disp{\sl\hfuzz 20pt\narrower\narrower\narrower\narrower
in each
bipartite graph, there exists a \emph{minimal vertex cover}
\(i.e., a vertex cover of the minimum possible cardinality\)
invariant with respect to all automorphisms of the graph
preserving the fractions.
}%
(Namely,
such a cover is the
set of all vertices of the fraction $A$ belonging to any minimal cover,
united with the set of vertices of the fraction $B$
belonging to at least one minimal cover.)

Let us choose such an invariant minimal vertex cover $Q'$
in $\Delta$ and note that each set $A_S$  
either is contained in $Q'$, or does not intersect $Q'$, 
because any permutation of~$A_S$ is an automorphism of the graph 
$\Delta$. So, 
$Q'=\(\bigsqcup\limits_{S\in\cal S'}A_S\)\sqcup(Q'\cap B)$
for some family~$\cal S'$.
For each $S\in\cal S'$, 
let us apply Theorem KL and choose a set $Z_S$
of edges
of the induced graph 
$\gp S$ 
\-
invariant with respect to all automorphisms 
of the graph $\gp S$, 
\-
representing all subgraphs (in $\gp S$) isomorphic $K_0$,
\- 
and such that $|Z_S|\le m_S\cdot|E(K_0)|$.

\enditem
The set $Y''=\(\bigcup\limits_{S\in\cal S'}Z_S\)\sqcup(Q'\cap B)$
is as required.
Indeed, 
it is, obviously, $(\Aut\Gamma')$-invariant,
intersects all subgraphs (in $\Gamma'$) isomorphic to $K$,
and 
$$
\eqalign{
|Y''|&=\sum_{S\in\cal S'}|Z_S|+|Q'\cap B|
\le
\sum_{S\in\cal S'}m_S\cdot|E(K_0)|+|Q'\cap B|
=
\sum_{S\in\cal S'}|A_S|\cdot|E(K_0)|+|Q'\cap B|
=
\cr
&=
\sum_{S\in\cal S'}|A_S|\cdot(|E(K)|-1)+|Q'\cap B|
\le
\(\sum_{S\in\cal S'}|A_S|+|Q'\cap B|\)\cdot(|E(K)|-1)
=
|Q'|\cdot(|E(K)|-1)
=
\cr
&=
|Q|\cdot(|E(K)|-1)
\le
|X'|\cdot(|E(K)|-1).
}
$$
By virtue of $(*)$, this completes the proof.

\s 3.
Proof of the theorem on symmetrisation of systems of 
multiple representatives

The reasoning below copies the argument from [KL21] almost literally;
the only difference is that, instead of B.~Neumann's theorem
[Neu54], we use its generalisation due to Sun (see below).

For $m=\max\limits_{F\in\F}|F|$, consider the following set
$
Y=\left\{y\in U\;\Bigm|\;
|(G\o y)\cap X| \ge{k\over m}|G\o y|\right\}
$
(in particular, $Y$ contains no points with infinite orbit).
Clearly, this set is $G$-invariant and
$k|Y|\le m|X|$ (because, for each orbit $G\o u$, we have
the inequality $k|(G\o u)\cap Y|\le m|(G\o u)\cap X|$).

It remains to show that $Y$ is a system of representatives for $\F$.
Take a set $F\in\F$.
Each set $g\o F$ (where $g\in G$) belongs to $\F$, since
$\F$ is invariant, and, hence, $|X\cap g\o F|\ge k$.
Therefore, the group $G$ is $k$-times covered by the family of sets
$$
G_f=\{g\in G\;|\;g\o f\in X\},
\qbox{where $f\in F$}.
$$
Each set $G_f$ is
either empty, or the union of finitely many left cosets
of the stabiliser $\St(f)$ of the point $f$:
$$
G_f=\{g\in G\;|\;g\o f\in X\}=
\bigcup_{x\in X}\{g\in G\;|\;g\o f=x\}=\!\!
\bigcup_{x\in X\cap G\o f}g_x\cdot \St(f),
\hbox{ where $g_x\in G$ are chosen such that $g_x\o f=x$}.
$$
Thus, we obtain a $k$-times cover of $G$
by
left cosets of some subgroups.
Let us
use a result of Sun ([Sun01], Lemma~2.2
+ [Sun90], Corollary 1%
\fn{%
In [Sun01], this was stated for finite-index subgroups;
while in [Sun90], it was shown that
cosets of
infinite-index subgroups may be
removed, and a $k$-times cover remains a $k$-times cover.
}%
):
\disp{\sl
if a group $G$ is
$k$-times
covered by finitely many cosets of its
\(not necessary different\) subgroups,
\newline
then
\
$\displaystyle\sum {1\over|G:G_i|}\ge k$
\(where the inverse to an infinite cardinal
is assumed to be zero\).
}%
Therefore,
(taking into account that the index of the stabiliser equals
the length of the orbit)
we obtain
$$
k\le\sum_{f\in F}{1\over|G:\St(f)|}\cdot|G\o f\cap X|=
\sum_{f\in F}{|G\o f\cap X|\over|G\o f|}.
$$
Since the number of terms in this sum is $|F|\le m$, some term must be at
least $k/m$, i.e., ${|G\o f\cap X|/|G\o f|}\ge k/m$, which means $f\in Y$
(by the definition of $Y$) and completes proof.

\s 4.
Proof of the theorem on
symmetrisation of systems of weighted representatives

We assume that no weight $w\in W$ exceeds one, because replacing each
weight function~$F\in\F$ with the 
function~
$x\mapsto\min\bigl(F(x),1\bigr)$
does not affect the notion of a 
system of weighted representatives.

Let us take an auxiliary finite set $E$ (for modelling weights) such that 
$|E|\cdot w$ is integer for all $w\in W$. The direct product 
$\~G=G\times S(E)$ of the group $G$ and the symmetric group~$S(E)$ of 
permutations of~$E$ acts naturally on the Cartesian product 
$\~U=U\times E$ 
(with projections $\pi_U\:\~U\to U$ and $\pi_E\:\~U\to E$).

Since the set $X\subseteq U$ is a system of weighted
representatives for the family of functions $\F$,
the set $\~X=X\times E$ is a system of representatives
of multiplicity $|E|$
for the following $\~G$-invariant
family of subsets of
$\~U$:
$$
\~\F
=
\bigl\{\~F\subseteq\~U\;\bigm|\;
\exists F\in\F\ \forall u\in U\
|\pi_U^{-1}(u)\cap\~F|=F(u)\cdot|E|
\bigr\}.
$$
In other words, for each weight function $F\:U\to\Q$,
we put into the set
$\~\F$ all subsets of the Cartesian product
$\~U=U\times E$ containing $F(u)|E|$ elements of $E$
``over" each point $u\in U$.
In the
student-council example (see Introduction),
we can take a
three-element set $E=\{1,2,3\}$; then the family
$\~\F$ consists of all subsets of the set
$\{students\}\times E$ of the form
$$
\bigl\{({\bf A},1),\;({\bf A},2),\;({\bf A},3),\;
({\bf B},i_1),\;({\bf C},i_2),\dots,\;({\bf J},i_9)\bigr\},
$$
where
$\bf B,C,D,E,F,G,H,I,J$ are nine different students respecting
another student $\bf A$ and $i_j\in\{1,2,3\}$.

By the theorem on the symmetrisation
of systems of multiple representatives, there exists a $\~G$-invariant
system of $|E|$-multiple
representatives $\~Y$ for this family such that
$$
|\~Y|
\le
|\~X\cap\~Y|\cdot\max\limits_{\~F\in\~\F}|\~F|/|E|
=
|\~X\cap\~Y|\cdot\max\limits_{F\in\F}\sum\limits_{u\in U} F(u).
\eqno{(3)}
$$

It remains to note that the projection
$Y=\pi_U(\~Y)\subseteq U$
is a required $G$-invariant system of weighted representatives for the
family of functions $\F$.
Indeed,
\-
we have
$
|Y|
=
|\~Y|/|E|
\le
|\~Y\cap\~X|\cdot\max\limits_{F\in\F}\sum\limits_{u\in U} F(u)/|E|
=
|Y\cap X|\cdot\max\limits_{F\in\F}\sum\limits_{u\in U} F(u)
$
(where the first equality follows from the $S(E)$-invariance
of $\~Y$;
the inequality is estimate (3);
and the last
equality $|\~Y\cap\~X|/|E|=|Y\cap X|$ is obvious,
because $\~X=X\times E$ (by the definition of $\~X$),
and $\~Y=Y\times E$ (by the definition of $Y$ and by the
$S(E)$-invariance of $\~Y$);
\-
$Y$ is a system of weighted representatives
for the family $\F$, because
$\~Y=Y\times E$ is a system of $|E|$-multiple representatives for the
family $\~\F$.

\References

[DM58]
A. L. Dulmage, N. S. Mendelsohn,
Coverings of bipartite graphs,
Canadian Journal of Mathematics, 10 (1958), 517-534.

[KhM07a]
E. I. Khukhro, N. Yu. Makarenko,
Large characteristic subgroups satisfying multilinear commutator
identities,
J. London Math. Soc., 75:3 (2007), 635-646.

[KhM07b]
E. I. Khukhro, N. Yu. Makarenko,
Characteristic nilpotent subgroups of bounded co-rank and
\newline
automorphically-invariant ideals of bounded codimension in Lie algebras,
{Quart. J. Math.}, 58 (2007), 229-247.

[KL21]
A. A. Klyachko, N. M. Luneva,
Invariant systems of representatives, or the cost of symmetry,
Discrete Mathematics, 344:6 (2021), 112361.
\arXiv 1908.03315

[KlMe09]
A. A. Klyachko, Yu. B. Mel'nikova,
A short proof of the Khukhro-Makarenko theorem on large characteristic
subgroups with laws,
Sbornik: Mathematics, 200:5 (2009), 661-664.
\arXiv 0805.2747

[KlMi15]
A. A. Klyachko, M. V. Milentyeva,
Large and symmetric: The Khukhro-Makarenko theorem on laws --- without laws,
J. Algebra, 424 (2015), 222-241.
\arXiv 1309.0571

[LP86]
L. Lov\'asz, M. D. Plummer,
Matching theory.
Elsevier, 1986.

[Ma71]
W. Mader,
Minimale $n$-fach kantenzusammenh\"angende Graphen,
Math. Ann. 191:1 (1971), 21-28.

[Neu54]
B. H. Neumann,
Groups covered by permutable subsets,
J. London Math. Soc., s1-29:2 (1954), 236-248.

[Sun90]
Zhi-Wei Sun,
Finite coverings of groups,
Fund. Math., 134:1 (1990), 37-53.

[Sun01]
Zhi-Wei Sun,
Exact $m$-covers of groups by cosets,
European Journal of Combinatorics, 22:3 (2001), 415-429.

[T22]
M. S. Terekhov,
The cost of symmetry in connected graphs,
Mathematical Notes, 112:6 (2022), 978-983.
\arXiv 2202.09590

\end